\documentclass[11pt]{article}

\begin{document}

\title{On a question of Kaplansky concerning his density theorem}
\author {George A. Elliott and Charles J. K. Griffin\\[4pt]
\small{\it Dedicated to Jacques Dixmier on the occasion of his 100th birthday} }

\date{}

\maketitle

\begin{abstract}
A new proof following a suggestion of Kaplansky to use a result of Dixmier, and in this way avoid unbounded nets of operators, is given of the Kaplansky density theorem.

\bigskip




\end{abstract}

\bigskip

In a remark concerning his fundamental density theorem (Theorem 1 of \cite{kap}), which states that the unit ball of a *-algebra of bounded Hilbert space operators is strongly dense in the unit ball of the strong closure of the algebra, Kaplansky reported (Remark 3(b) of \cite{kap}) that while his theorem immediately implied Dixmier's theorem that a *-algebra of operators is strongly closed if its unit ball is strongly closed (Th\'eor\`eme 8 of \cite{dix}), it did not appear to be possible to reverse the reasoning.  In fact, Kaplansky's insight was correct---the reverse implication holds.

This is seen as follows.  Given a *-algebra of operators on a Hilbert space---which we may assume to be norm-closed---instead of considering explicitly the strong operator closure of this algebra---the limits of all strongly convergent nets in the algebra---, consider just the limits of bounded nets.  It is easy to check (see below) that these form a *-algebra, and also (using only strong continuity of the continuous functional calculus on a bounded set of self-adjoint operators) that the unit ball of this *-algebra is the strong closure of the unit ball of the given *-algebra and in particular is strongly closed.  Hence by the implication $(f) \Rightarrow (b)$ of  Th\'eor\`eme 8 of \cite{dix}, the whole *-algebra is strongly closed, and in particular is the strong closure of the given algebra.  Almost incidentally to the proof that the unit ball of this algebra is strongly closed (see below), its unit ball is the strong closure of the unit ball of the given *-algebra---exactly the assertion of the Kaplansky density theorem.

To see that the strong limits of bounded nets in the given *-algebra form a *-algebra, recall that the operations of addition and scalar multiplication are jointly strongly continuous (without restriction), and the operation of multiplication is jointly strongly continuous on a  bounded set. As far as the *-operation is concerned, recall that strong and weak operator closures coincide for convex sets, and that the *-operation is weakly continuous. (All of these facts are also used in \cite{kap}.)

Finally, note that to prove that the unit ball of the larger *-algebra is strongly closed, asserted above, it is enough to prove that it is the strong closure of the unit ball of the given algebra---this fact is of course pertinent in itself! First, given a self-adjoint element of this unit ball, by definition a limit of a bounded net in the given *-algebra, by convexity and weak continuity of the *-operation as above, one expresses this element as the limit of a bounded net of self-adjoint elements of the given (norm-closed) *-algebra.  Acting on this net with a suitable real-valued continuous function (as in \cite{kap}, but now the net is bounded!) yields a net in the unit ball with the same strong limit.  The case of a non-self-adjoint element of the unit ball of the larger *-algebra is then reduced to the self-adjoint case by passing to $2 \times 2$ matrices as in \cite{kap} (a device attributed by Kaplansky to Halmos).

\medskip

\noindent{\bf Remark}
The present proof may be broken down into two ingredients---first, Dixmier's result (Th\'eor\`eme 2 of \cite{dix}) that the algebra of bounded operators on a Hilbert space (hence any von Neumann algebra) is in a natural way (Sakai later, in \cite{sak}, proved uniquely) a Banach space dual, with the weak$^*$ topology coinciding with the ultraweak operator topology, and second, the result of Dieudonn\'e (Th\'eor\`eme 23 of \cite{dix})---attributed by him to Bourbaki (\cite{bou}) and, in equivalent form, to Banach (\cite{ban})---that a linear subspace of a Banach space dual is weak$^*$ closed if its unit ball is. (The second result was generalized to convex sets in \cite{kre}.)

\smallskip

\noindent{\small Department of Mathematics, University of Toronto, Toronto, Ontario, Canada~~ M5S 2E4\\
e-mail: elliott@math.toronto.edu}

\smallskip

\noindent{\small Department of Mathematics, University of Toronto, Toronto, Ontario, Canada~~ M5S 2E4\\
e-mail: charlie.griffin@mail.utoronto.ca}

\end{document}